\documentclass[a4paper,oneside]{amsart}

\usepackage[textwidth=125mm,textheight=195mm]{geometry}
\usepackage{amssymb,amsthm,amsmath}
\usepackage{mathtools}
\usepackage{microtype}
\usepackage{graphicx}
\usepackage{float}
\usepackage{enumitem}
\usepackage{mathrsfs}
\usepackage{hyperref}

\theoremstyle{plain}
\newtheorem{thm}{Theorem}[section]
\newcommand{\ceqq}{\;=\;}
\newcommand{\cplus}{\;+\;}

\newcommand{\ccoloneqq}{\;\coloneqq\;}
\newcommand{\act}{\!\cdot\!}
\newcommand{\multiset}{\mathcal{M}}

\newcommand{\mc}[1]{\mathcal #1}
\newcommand{\TO}{\mc L}

\newcommand{\eps}{\varepsilon}
\newcommand{\wh}{\widehat}

\newcommand\N{\mathbb{N}}
\newcommand\Q{\mathbb{Q}}
\newcommand\R{\mathbb{R}}
\newcommand\Z{\mathbb{Z}}
\newcommand\C{\mathbb{C}}
\newcommand\T{\mathbb{T}}
\newcommand{\h}{\mathbb{H}}
\newcommand{\mP}{\mathbb{P}}

\newcommand{\parab}{\textnormal{par}}
\newcommand{\Vect}{\mc V_s^{\star}}

\newcommand{\mat}[4]{\begin{pmatrix} #1&#2\\#3&#4\end{pmatrix}}
\newcommand{\bmat}[4]{\begin{bmatrix} #1&#2\\#3&#4\end{bmatrix}}
\newcommand{\textmat}[4]{\left(\begin{smallmatrix} #1&#2 \\ #3&#4
\end{smallmatrix}\right)}
\newcommand{\textbmat}[4]{\left[\begin{smallmatrix} #1&#2 \\ #3&#4
\end{smallmatrix}\right]}

\DeclareMathOperator{\id}{id}
\DeclareMathOperator{\SL}{SL}
\DeclareMathOperator{\PSL}{PSL}
\DeclareMathOperator{\Ima}{Im}
\DeclareMathOperator{\Rea}{Re}
\DeclareMathOperator{\tr}{tr}
\DeclareMathOperator{\acosh}{acosh}
\DeclareMathOperator{\HP}{HP}
\DeclareMathOperator{\base}{base}

\begin{document}
\title[Dynamics and Maass forms]{Dynamics of geodesics, and Maass cusp forms}

\author{Anke Pohl} 
\address{Anke Pohl, University of Bremen, Department~3 -- Mathematics, Bibliothekstr.~5, 28359 Bremen, Germany}
\email{apohl@uni-bremen.de}

\author{Don Zagier}
\address{Don Zagier, Max Planck Institute for Mathematics, Vivatsgasse~7, 53111 Bonn, Germany
and
International Centre for Theoretical Physics, Strada Costiera, Trieste, Italy}
\email{dbz@mpim-bonn.mpg.de}

\begin{abstract}
The correspondence principle in physics between quantum mechanics and classical mechanics suggests deep relations between spectral and geometric entities of Riemannian manifolds. We survey---in a way intended to be accessible to a wide audience of mathematicians---a mathematically rigorous instance of such a relation that emerged in recent years, showing a dynamical interpretation of certain Laplace eigenfunctions of hyperbolic surfaces.  
\end{abstract}

\maketitle

\section{Introduction}

Suppose we have a huge space, such as the earth or a billiard table, and a small marble sitting on this space. We give this marble an initial push and observe its trajectory as it travels over the space. As we experienced from a very young age on, the marble goes straight until it hits an obstacle, e.\,g., the boundary of the billiard table, from which it bounces off with outgoing angle equal to incoming angle, and then continues its straight path until the next obstacle where the same game restarts.

\begin{figure}[h]
\centering
\includegraphics{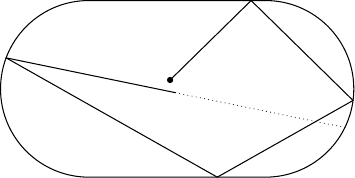} 
\caption{Trajectory on a stadium-shaped billiard table.}\label{fig:billiard}
\end{figure}

\begin{figure}
\centering
\includegraphics{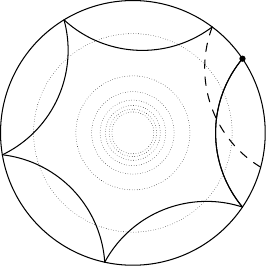} 
\caption{Trajectory on a hill, given by a disk with a bump in the middle, viewed from above. Height level curves are indicated by dotted circles.}\label{fig:bump}
\end{figure}

In Figure~\ref{fig:billiard} this situation is depicted for a flat stadium-shaped billiard table. In Figure~\ref{fig:bump} it is shown for a disk with a bump in the middle, indicating that `straight path' here means `path of minimal resistance' or `path of minimal effort'.

In terms of physics, the motion of the marble is predicted by the laws of classical mechanics. In such a description, moving objects
are often modeled as point particles, that is, as objects without size or dimension, identifying the object with its center of mass.

In reality, any real-world object has a non-zero size, and the idealization as a point is not always desirable or correct. If we consider a very small marble which is almost a point, say of the size of an electron, or if we zoom in into our previous marble and try to describe the trajectory of a single electron of it then we notice that the classical mechanics model is not accurate on this subatomic level. One of the obstacles is the impossibility to determine simultaneously with absolute precision the position and momentum of the considered particle, as expressed by Heisenberg's famous uncertainty principle. Thus, the classical mechanical principles of determinism and time reversibility are not valid anymore. On such small scale, a more accurate model is provided by quantum mechanics, which describes the probability with which the particle attains a specific position-momentum combination.

The correspondence principle in physics states that, in the limit of passing to large scale, the predictions of quantum mechanics reproduce those of classical mechanics. However, the precise relation between classical and quantum mechanics is not yet fully understood, and its investigation gives rise to many interesting mathematical questions. 

In terms of mathematics, the classical mechanical aspects of the motion of the marble considered above translate to properties of the geodesic flow on a Riemannian manifold~$X$, whereas the quantum mechanical description relates to the Laplace operator on~$X$ and its ($L^2$-)eigenvalues and eigenfunctions. The correspondence principle then suggests an intimate relation between geometric-dynamical aspects of~$X$ on the one hand, and its spectral aspects on the other:

\begin{center}
 \begin{tabular}{rcrl}
 \multicolumn{1}{c}{\textbf{physics}} &  & \multicolumn{2}{c}{\textbf{mathematics}\hphantom{xxxxxxxxx}} 
  \\[1mm]\hline\\[-2mm]
  classical mechanics & $\leftrightsquigarrow$ &  geometric entities\,: \hspace*{-3mm}&$\begin{cases} \text{periodic geodesics} \\ \text{lengths of periodic geodesics} \end{cases}$
  \\[5mm]
  quantum mechanics& $\leftrightsquigarrow$ & spectral entities\,: & $\begin{cases} \text{Laplace eigenfunctions} \\ \text{Laplace eigenvalues} \end{cases}$
 \end{tabular}
\end{center}

During the last century, many results showing relations between geometric-dynamical and spectral properties of Riemannian manifolds have been obtained. In Section~\ref{sec:appetizer} we will discuss---as an appetizer---the flat $1$-torus where a clear relation between the lengths of periodic geodesics (`classical mechanical objects') and the Laplace eigenvalues (`quantum mechanical objects') appears. 

The main aim of this article is to present a much deeper relation between periodic geodesics and Laplace eigenfunctions that has emerged in recent years, but now for a class of hyperbolic surfaces.

In a nutshell, this goes as follows. A well-chosen discretization of the flow along the periodic geodesics gives rise to a one-parameter family of \emph{transfer operators}, which are evolution operators that are reminiscent of weighted graph Laplacians and that also may be thought of as discretizations of the hyperbolic Laplacian. As such, these operators are simultaneously objects of classical and quantum mechanical nature, and therefore can serve as mediators between the dynamical and spectral entities of the hyperbolic surface under consideration. In our case, highly regular, rapidly decaying eigenfunctions (called \emph{period functions}) of eigenvalue~$1$ of the transfer operator with parameter~$s$ are in bijection with rapidly decaying Laplace eigenfunctions  (called \emph{Maass cusp forms}) with spectral parameter~$s$. This provides a purely dynamical interpretation of the Maass cusp forms (not just their eigenvalues), shows a close dependence between periodic geodesics and these Laplace eigenfunctions, and provides a deep-lying mathematical realization of an instance of the correspondence principle. 

The modular surface was the first hyperbolic surface for which such a result could be established, through combination of work by E.~Artin~\cite{Artin}, Series~\cite{Series}, Mayer~\cite{Mayer_thermo, Mayer_thermoPSL}, Lewis~\cite{Lewis}, Bruggeman~\cite{Bruggeman_lewiseq}, Chang--Mayer~\cite{Chang_Mayer_transop}, and Lewis--Zagier~\cite{LZ_survey, LZ01}.
Taking advantage of the constructions involved, an extension to a class of finite covers of the modular surface was achieved in the combination of~\cite{Chang_Mayer_transop, Deitmar_Hilgert, Fraczek_Mayer_Muehlenbruch}. An alternative proof for the modular surface was provided in~\cite{BM09, Mayer_Muehlenbruch_Stroemberg}. The recent development of a new type of discretizations for geodesic flows on hyperbolic surfaces~\cite{Pohl_Symdyn2d} and of a cohomological interpretation of the Maass cusp forms~\cite{BLZm} allowed to prove such a relation between periodic geodesics and Laplace eigenfunctions for a large class of hyperbolic surfaces far beyond the modular surface and in a very direct way~\cite{Moeller_Pohl, Pohl_mcf_general, Pohl_mcf_Gamma0p}. 

In Sections~\ref{sec:mod_surface}--\ref{sec:recap} we will survey this new approach, although in an informal way and restricting for simplicity to the modular surface. We attempt to provide sufficiently precise definitions and enough details to keep the exposition as understandable as possible without introducing too much technical material. As a general principle we invite all readers to rely on their intuitive understanding of the geometry and dynamics of Riemannian manifolds, to use the many figures as a support, and to ignore the exact expressions of all formulas.

To end this introduction, we briefly mention another example of the many other strands of research seeking for and establishing relations between geometric-dynamical and spectral properties of Riemannian manifolds, where considerable progress has been made in the last two decades and that provides another concrete incarnation of the correspondence principle: the problem of quantum unique ergodicity. This problem concerns the distribution of the mass of high energy Laplace eigenfunction (i.e., with large eigenvalue). A conjecture by Rudnick and Sarnak states that on surfaces with sufficiently chaotic geodesic flows, the mass of Laplace eigenfunctions equidistributes as their eigenvalues tend to infinity. In other words, for such surfaces, the limiting behavior of the mass distribution of Laplace eigenfunctions is expected to be governed by the behavior of the geodesic flow. We refer to~\cite{Hassell, Sarnak, Zelditch} for precise statements and excellent surveys of the recent developments.

\noindent
\emph{Acknowledgements.} 
AP wishes to thank the Max Planck Institute for Mathematics in Bonn for hospitality and excellent working conditions during the preparation of this manuscript. Further, she acknowledges support by the DFG grants PO~1483/2-1 and PO~1483/2-2.

\section{An appetizer}\label{sec:appetizer}

In this section we will treat the `baby case' of the flat \emph{$1$-torus}
\[
 \T \ceqq \R/\Z \ceqq [0,1]/{\{0\!=\!1\}}\,,
\]
and show an intimate and very clear relation between geometric and spectral entities, and hence a mathematical rigorous instance of the correspondence principle.

Of course, this specific one-dimensional Riemannian manifold is much too simple to be representative of the general situation. However, it allows us to provide---without too much technical effort---a first instance of the relation between the geometry and the spectrum as motivated by the considerations from physics. We will also use this `baby example' to carefully introduce the relevant geometrical and spectral concepts, whose counterparts in the situation of hyperbolic surfaces will be treated in the main body of this paper.

\subsection{The flat 1-torus.}

\begin{figure}[h]
\centering
\includegraphics{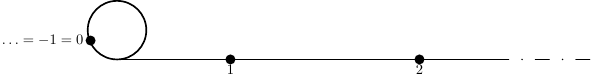} 
\caption{Rolling up $\R$ to form $\T$.}\label{fig:torus}
\end{figure}

For a pictorial, but rather sketchy construction of the flat $1$-torus~$\T$ we may imagine the set~$\R$ of real numbers as a number line, and glue together this line at any two points that are separated by an integer distance. The glueing process can be visualized as rolling up the line to a unit circle. (See Figure~\ref{fig:torus}.) Alternatively, we may take the interval~$[0,1]$ and glue together its two endpoints~$0$~and~$1$. 

Both these geometric constructions indicate that $\T$ carries more structure than just being a set. In particular, as we will explain now, a well-defined notion of distances on~$\T$ exists, and derivatives of maps from and to~$\T$ can be defined. 

In order to be able to formulate such additional structures in precise terms and to work with them, we use a formula-based definition of~$\T$. For that, we identify any two points of~$\R$ that differ by an integer only. Thus, for each~$r\in\R$, all points in the set
\begin{equation}\label{def:coset}
 \{ r+m \mid m\in\Z\}
\end{equation}
are unified to a single element, which we denote by~$[r]$. The torus~$\T$, as a set, consists of all these elements. The glueing process in the pictorial construction is a visualization of the \emph{projection map}
\begin{equation}\label{eq:projmap}
 \pi_{\T }\colon \R\to \T \,,\quad r\mapsto [r]\,.
\end{equation}
This map is \emph{locally injective}, which means that for any~$r\in\R$ we find a small~$\eps>0$ such that the restriction of~$\pi_\T$ to the interval~$(r-\eps,r+\eps)$ is injective. Here, we may choose $\eps=\tfrac12$ for each~$r\in\R$. In rough terms, small pieces of the torus~$\T$ look exactly like small pieces of~$\R$. It is precisely this property which allows us to push certain structures of~$\R$ to~$\T$.

\subsection{Geometric entities.}\label{sec:geom_entities}

The \emph{geometric entity} or, from the standpoint of the introduction, the \emph{classical mechanical object} that we are interested in is the set of all periodic geodesics. A geodesic is a path on~$\T$ of a specific type that we now introduce. 

In order to define the differentiability of a function mapping from an open interval in~$\R$ to~$\T$ we pick for each~$x\in\T$ a representative~$r_x\in\R$ (thus, $[r_x] = x$) and denote by~$\sigma_x$ the inverse of the restriction of the projection map~$\pi_\T$ to the interval~$(r_x-\tfrac12,r_x+\tfrac12)$. (We recall that $\pi_\T$ is locally injective.) Then $\sigma_x$ is the bijective map 
\[
 \sigma_x\colon \T\smallsetminus\{ [r_x+\tfrac12]\} \to \bigl(r_x-\tfrac12,r_x+\tfrac12\bigr)
\]
satisfying
\[
 \pi_\T\circ\sigma_x = \id_{\T\smallsetminus\{ [r_x+1/2]\}}\,.
\]
Let $I\subseteq\R$ be an open interval and $p\colon I\to\T$ a map. Then~$p$ is \emph{differentiable} at~$t\in I$ if there exists~$\delta>0$ such that the map 
\[
 \sigma_{p(t)}\circ p\colon (t-\delta, t+\delta) \to \R
\]
is well-defined and differentiable at~$t$. In this case, the \emph{derivative} of~$p$ at~$t$ is 
\[
 p'(t) \ccoloneqq \bigl(\sigma_{p(t)}\circ p\bigr)'(t)\,.
\]
It is straightforward to check that neither the property of differentiability nor the derivative depends on the choice of the representative of~$p(t)$ in~$\R$. The map~$p$ is a \emph{path} on~$\T$ if it is differentiable at any point of~$I$, a property also called \emph{differentiable} for short. For any path~$p\colon I\to\T$, the set $I$ should be thought of as a time interval, and $p(t)$ as the position where we are at time $t$ if we travel along the path $p$. The derivative~$p'$ is the \emph{speed} of~$p$, and $p$ is said to be of \emph{unit speed} if $|p'(t)|=1$ for all~$t\in I$. 

A path~$p\colon I\to\T$ is \emph{straight} or a \emph{geodesic} on~$\T$ if---roughly said---for any two nearby points on the path no shorter way between them exists than the path itself. To be more precise, we define the \emph{distance} between two points~$x,y\in\T$ to be the minimal distance between any two of their representatives in~$\R$, hence
\[
 d_\T(x,y) \ccoloneqq \min\big\{ d_\R(r_x,r_y) \ \big\vert\ [r_x] = x,\ [r_y]=y\big\}\,,
\]
where 
\[
 d_\R(r_x,r_y)  \ccoloneqq | r_x -r_y |
\]
is the usual euclidean distance on~$\R$. A path~$p\colon I\to\T$ of unit speed is \emph{straight} if for any~$t\in I $ there exists~$\eps>0$ such that for all~$t_1,t_2\in (t-\eps,t+\eps)\cap I$ we have
\[
 d_\T\big(p(t_1),p(t_2)\big) \ceqq \left| \int_{t_1}^{t_2} |p'(t)|\, dt\right| \ceqq \big|t_1-t_2\big| \ceqq d_\R(t_1,t_2)\,.
\]
That is, the distance between~$p(t_1)$ and~$p(t_2)$ equals the length of the path between~$p(t_1)$ and~$p(t_2)$, which here also equals the euclidean distance between~$t_1$~and~$t_2$. From now on, `geodesic' will always mean a \emph{unit speed, complete geodesic}, i.\,e., a straight path of unit speed with time interval~$I=\R$. 

In everyday language, the notion of path usually does not refer to the motion, i.\,e., to a map $p\colon I \to\T$, but rather to the static object, i.\,e., to the image~$p(I)$ of~$p$. The orientation, however, is important: `the path from~$a$ to~$b$'. We too will use the notion of geodesic more flexibly and apply it to refer to either
\vspace*{-.5mm}
\begin{enumerate}[label=$\mathrm{(\mathbf{G\arabic*})}$, ref=$\mathrm{(\mathbf{G\arabic*})}$]
 \item\label{geod1} a geodesic~$p\colon\R\to\T$ defined as above as a path, or
 \item\label{geod2} \parbox[t]{\textwidth-1.5cm}{the oriented image of such a geodesic, or---more precisely---its equivalence class when we identify any two such geodesics that differ only by a shift in their arguments.}
\end{enumerate}
The motivation for the second usage is that we are typically not interested in the specific time parametrization of a geodesic. The context should always clarify which version is being used. 

\begin{figure}[h]
\centering
\includegraphics{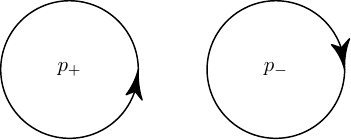} 
\caption{The two periodic geodesics on $\T$.}\label{fig:geod_torus}
\end{figure}

In our one-dimensional `baby example' there are only two geodesics in the sense~\ref{geod2}, namely those represented by the two geodesics in the sense~\ref{geod1} given by
\[
 p_\pm \colon \R\to\T ,\quad t\mapsto [\pm t]\,.
\]
(See Figure~\ref{fig:geod_torus}.) 
Both these geodesics are \emph{periodic}, that is, they `close up', or in rigorous terms, there exists~$t_0 > 0$ such that for all~$t\in \R$, 
\[
 p_\pm(t) \ceqq p_\pm (t+t_0)\,.
\]
The minimal such~$t_0$ is called the \emph{(primitive) period} or \emph{(primitive) length}~$\ell(p_\pm)$ of the geodesic~$p_\pm$, which here is $\ell(p_\pm)=1$ in both cases. Periodicity and lengths are invariants under the equivalence of geodesics, and hence an intrinsic notion for geodesics in the sense~\ref{geod2}. 

For this `baby example' we are interested in the \emph{(primitive) geodesic length spectrum}~$L_\T$, defined as the multiset ($=$ set with multiplicities) of lengths of the periodic geodesics in the sense~\ref{geod2}. In our case, this is
\[
 L_{\T } \ceqq\{\text{lengths of periodic geodesics}\}_\multiset \ceqq \{1,1\}_\multiset\,.
\]

\subsection{Spectral entities.}

The \emph{spectral entity} or the \emph{quantum mechanical object} that we need here is the Laplace spectrum of~$\T $, which we now explain.

In analogy with the definitions of differentiability for maps~$I\to\T$,  we use the local injectivity of the projection map~$\pi_\T$ to characterize the differentiability of functions~$\T\to\C$. Then a function $f\colon \T \to \C$ is \emph{differentiable} if the map
\[
 F\ccoloneqq f\circ\pi_\T \colon \R \to \C
\]
is differentiable. The \emph{derivative} of~$f$ at~$[r]\in\T$ is then the derivative of~$F$ at~$r\in\R$, and again it is straightforward to check that $f'([r])$ is indeed well-defined.

Further, a function $f\colon\T\to\C$ is \emph{square-integrable} or, for short, in~$L^2$, if the map
\[
 \tilde f\ccoloneqq f\circ\pi_\T\vert_{[0,1)} \colon [0,1) \to\C
\]
is square-integrable in the usual sense. In particular, $\tilde f$ is integrable and 
\[
 \int_0^1 |\tilde f(r)|^2\,dr < \infty\,.
\]
We identify any two $L^2$-functions $f_1,f_2\colon\T\to\C$ for which $\int_0^1|\tilde f_1(r)-\tilde f_2(r)|^2\,dr = 0$, and denote the set of all equivalence classes by~$L^2(\T)$.

The \emph{Laplace operator} on~$\T $, given by
\[
 \Delta_{\T } \ccoloneqq -\frac{d^2}{d [r]^2}\,,
\]
acts on~$L^2(\T)$. As Fourier theory shows, a basis for its $L^2$-eigenfunctions is constituted by the family
\[
 f_k\colon \T  \to \C\,,\quad f_k\big([r]\big) \coloneqq e^{2\pi i k r}\qquad (k\in\Z)\,.
\]
An immediate calculation gives
\[
 \Delta_{\T } f_k \ceqq (2\pi k)^2 f_k\,.
\]
Thus, the \emph{Laplace spectrum} of~$\T $ is the multiset
\[
 \sigma(\T )\ceqq \{\text{Laplace eigenvalues}\}_\multiset \ceqq \{ (2\pi k)^2 \mid k\in\Z\}_\multiset\,.
\]

\subsection{Relation between geometric and spectral entities.}

The physics-informed intuition on a close relation between the geodesic length spectrum~$L_\T$ of~$\T$ and the Laplace spectrum~$\sigma(\T)$ can be proven mathematically rigorously in different ways of which we provide one here. For that we consider the \emph{dynamical zeta function}
\[
 \zeta_{\T }(s) \ccoloneqq \prod_{\ell\in L_{\T }}\left( 1 - e^{-s\ell}\right) \ceqq \left( 1 - e^{-s}\right)^2\,.
\]
Then
\[
 \zeta_{\T }(s) = 0 \qquad\Longleftrightarrow\qquad s=2\pi i k \quad\text{for some $k\in\Z$\,,}
\]
and the order of each zero is $2$. In other words,
\begin{equation}\label{eq:relation}
 \zeta_{\T }(s) = 0 \qquad\Longleftrightarrow\qquad (is)^2 \in \sigma(\T )\,,
\end{equation}
and the order of $s$ as a zero corresponds to the order of~$(is)^2$ as eigenvalue, except for $s=0$, where the order of the Laplace eigenvalue $(is)^2 = 0$ is $1$, whereas the order of the zero~$s=0$ of~$\zeta_\T$ is~$2$. 

Thus knowing the geodesic length spectrum~$L_\T$, and hence the dynamical zeta function~$\zeta_\T$, we can deduce all Laplace eigenvalues, and even their multiplicities up to the difficulty at~$s=0$. Conversely, if we are given the Laplace spectrum~$\sigma(\T)$ (with multiplicities), and hence all zeros of~$\zeta_\T$ with almost all multiplicities, then we can easily deduce the exact formula of~$\zeta_\T$ and thus the geodesic length spectrum. 

This ends the $1$-dimensional `appetizer'. In the rest of the paper we will study a $2$-dimensional case, again describing first the geometric side, then the spectral side, and then the relation between them. Of course, this case is much more involved, but we have tried to introduce the concepts in this one-dimensional torus case in such a way that they generalize naturally. 

\section{Geometric and spectral sides of the modular surface}\label{sec:mod_surface}

In the previous section we considered the torus~$\T$, which is a quotient of the flat \mbox{$1$-}manifold~$\R$ by a discrete group action. From now on, we will consider hyperbolic surfaces, which are orbit spaces of the \emph{hyperbolic plane}~$\h$ by discrete groups of Riemannian isometries. For concreteness we will discuss only the \emph{modular surface}~$X=\PSL_2(\Z)\backslash\h$, even though the results hold for a much larger class. We will provide precise definitions of all objects further below in this section.

In the course of the following four sections we will survey---as already mentioned in the introduction---a rather deep relation between the geodesic flow on~$X$ and the rapidly decaying Laplace \mbox{$L^2$-}eigenfunctions for the modular group~$\PSL_2(\Z)$, the \emph{Maass cusp forms}. This results in a dynamical interpretation of Maass cusp forms, or from a physics point of view, in a description of certain quantum mechanical wave functions using only tools and objects from classical mechanics. The proof of this relation is split into three major steps:
\begin{enumerate}[label=$\mathrm{(\Roman*)}$, ref=$\mathrm{\Roman*}$]
 \item A cohomological interpretation of Maass cusp forms, which we will explain in Section~\ref{sec:mcf}. Representing Maass cusp forms faithfully as cocycle classes in suitable cohomology spaces provides an interpretation of these forms in a rather algebraic way of which we will take advantage. 
 \item A well-chosen discretization of the geodesic flow on~$X$, which we will construct in Section~\ref{sec:discretization}. This discretization extracts those geometric and dynamical properties from the geodesic flow on~$X$ that are crucial for the relation to Maass cusp forms, and it discards all the other additional properties. This condensed, discrete version of the geodesic flow is also of a rather algebraic nature.
 \item\label{step3} A connection between the discretization of the geodesic flow and the cohomology spaces, as discussed in Section~\ref{sec:TO}. The central object mediating between these objects is the evolution operator (with specific weights, adapted to the spectral parameter of Maass cusp forms; a \emph{transfer operator}) of the action map in the discrete version of the geodesic flow. We will see that the highly regular eigenfunctions of the evolution operator with parameter~$s$ are building blocks for the cocycle classes in the cohomological interpretation of the Maass cusp forms with spectral parameter~$s$, and will establish an explicit bijection between these eigenfunctions and the Maass cusp forms.
\end{enumerate}

Even though the first two steps are technically independent of each other, crucial choices in the construction of the discretization of the geodesic flow in the second step can be motivated by the precise expressions in the cohomological interpretation of Maass cusp forms in the first step. Therefore we recommend the reader to go through these steps in the order as presented. The third step necessarily takes advantage of the results from Sections~\ref{sec:mcf} and \ref{sec:discretization}. From a technical point of view, only the final results of these sections are needed for Step~\eqref{step3}, not the information on how they were obtained, so readers who are only interested in this step may proceed directly to Section~\ref{sec:TO} after familiarizing themselves with the general setup and Theorems~\ref{thm:BLZ} and~\ref{thm:CS}. In Section~\ref{sec:recap} we will provide a brief recapitulation.

In the remainder of this section we introduce the geometric and spectral objects that we will need further on. We restrict ourselves here to the absolutely necessary minimum. There are many excellent textbooks which provide much more detail on these objects and comprehensive treatments of hyperbolic surfaces. We refer in particular to~\cite{Bergeron, Ratcliffe, Venkov}.

\subsection{The hyperbolic plane.}
The \emph{hyperbolic plane} is a certain two-dimensional manifold with Riemannian metric in which Euclid's parallel axiom fails: on the hyperbolic plane, for every straight line~$L$ (infinitely extended in both directions) and any point~$p$ not on~$L$ there are infinitely many lines~$\tilde L$ passing through~$p$ that do not intersect~$L$.

Abstractly, the hyperbolic plane is the unique two-dimensional connected, simply connected, complete Riemannian manifold with constant sectional curvature~$-1$ (see, e.g.,~\cite[Theorem~6.3]{Boothby}). There are many models for the hyperbolic plane. We use its \emph{upper half plane model}\footnote{Another widely known model for the hyperbolic plane is the Poincar\'e disk model, which prominently features in several of M.\@ C.\@ Escher's pictures.}
\[
 \h \ccoloneqq \{z\in\C \mid \Ima z > 0\}\,,
\]
where the \emph{line element of the Riemannian metric} is given by 
\begin{equation}\label{eq:lineelement}
ds^2_{x+iy} \ccoloneqq \frac{dx^2 + dy^2}{y^2}\,.
\end{equation}
Informally, the Riemannian metric allows us to measure distances and angles. Angles in hyperbolic geometry are identical to the euclidean angles in~$\h$. Distances between points however are changed in hyperbolic geometry when compared to euclidean geometry. From a euclidean point of view, hyperbolic distances between two points increase when these move nearer to the real axis~$\R$.

\begin{figure}[h]
\centering
\includegraphics{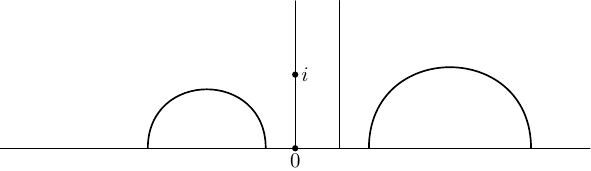} 
\caption{Geodesics on~$\h$.}\label{fig:hypplane}
\end{figure}

For the torus~$\T$ we discussed two notions of geodesics in Section~\ref{sec:geom_entities}: the \ref{geod1}-version in which we understand geodesics as paths, and the \ref{geod2}-version where we understand geodesics as oriented subsets. 
In the upper half plane model of the hyperbolic plane, the \ref{geod2}-version of geodesics, i.\,e., infinite paths that are straight with respect to this metric, are the (oriented) semi-circles with center on~$\R$ and the vertical rays based on the real axis. (See Figure~\ref{fig:hypplane}.) 

The upper half plane~$\h$ has a boundary whose definition is motivated by the dynamics of the geodesics on~$\h$; it consists of all `infinite endpoints' of the geodesics. Considering Figure~\ref{fig:hypplane}, this boundary is given by
\[
\mP^1(\R) \ccoloneqq \R\cup \{\infty\}\,.
\] 

A \emph{Riemannian isometry} is a bijective map on~$\h$ which preserves the distance between any two points. In particular, any Riemannian isometry maps geodesics to geodesics. 
The group of \emph{orientation-preserving Riemannian isometries} on the hyperbolic plane is isomorphic to the (projective) matrix group
\[
 G\ccoloneqq\PSL_2(\R) \ccoloneqq \SL_2(\R)/\{\pm \id\}\,.
\]
The element $g\in G$ represented by the matrix~$\textmat{a}{b}{c}{d}\in \SL_2(\R)$ is denoted by~$g =\textbmat{a}{b}{c}{d}$, with square brackets. It then has one other representative in~$\SL_2(\R)$, namely $\textmat{-a}{-b}{-c}{-d}$.
The \emph{action} of~$G$ on~$\h$ is given by
\begin{equation}\label{eq:action}
 \bmat{a}{b}{c}{d}\act z \ccoloneqq \frac{az+b}{cz+d}\,.
\end{equation}
Occasionally, we will omit the dot~$\cdot$ in the notation. The action of~$G$ on~$\h$, as defined in~\eqref{eq:action}, extends continuously to an action of~$G$ on~$\h\cup\nobreak \mP^1(\R)$ in the obvious way, using that in hyperbolic geometry the equality $1/0=\infty$ is valid. Thus, the right hand side of~\eqref{eq:action} is replaced by~$a/c$ if~$z=\infty$ and~$c\not=0$, and by~$\infty$ if $z=\infty$ and $c=0$ or if $cz+d=0$. We use the notation from~\eqref{eq:action} also for this extended action.

\subsection{The modular surface.}
\begin{figure}[h]
\centering
\includegraphics{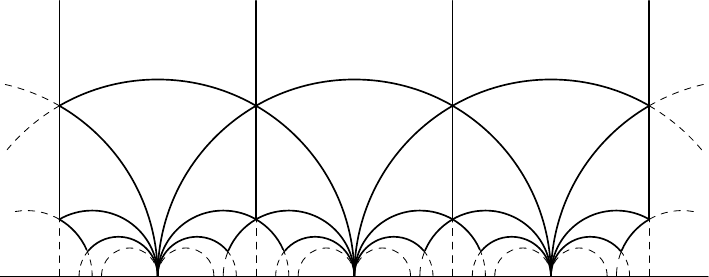} 
\caption{Tesselation of~$\h$ by triangles.}\label{fig:tess}
\end{figure}

A subgroup of~$G$ of particular importance is the \emph{modular group} 
\[
 \Gamma \ccoloneqq \PSL_2(\Z)\,.
\]
It acts on~$\h$ preserving the tesselation by triangles as indicated in Figure~\ref{fig:tess}.
\noindent
The \emph{modular surface} is the orbit space
\[
 X \ccoloneqq \Gamma\backslash\h\,,
\]
that is, the space we obtain if we identify any two points of~$\h$ that are mapped to each other by some element of~$\Gamma$. We let 
\begin{equation}\label{def:projmap}
 \pi\colon \h \to X=\Gamma\backslash\h
\end{equation}
be the projection map. The space~$X$ can be compactified by adding an additional point that is represented in~$\h$ by~$\infty~  (=i\infty)$. For future purpose we note that $\mP^1(\Q) \coloneqq \Q\cup\{\infty\}$ is the $\Gamma$-orbit of~$\infty$ and that the map~$\pi$ extends canonically to a map
\[
 \h\cup \mP^1(\Q) \to \overline X = \Gamma\backslash(\h\cup \mP^1(\Q))\,,
\]
which we continue to denote~$\pi$. 

\begin{figure}[h]
\centering
\includegraphics{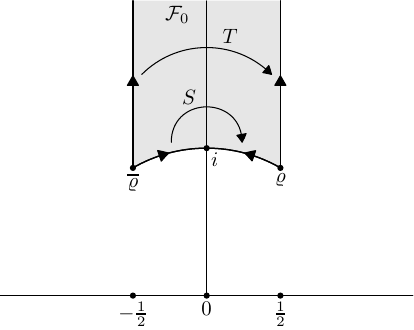} 
\caption{Fundamental domain~$\mc F_0$ for~$\Gamma$.}\label{fig:funddom1}
\end{figure}

\begin{figure}[h]
\centering
\includegraphics{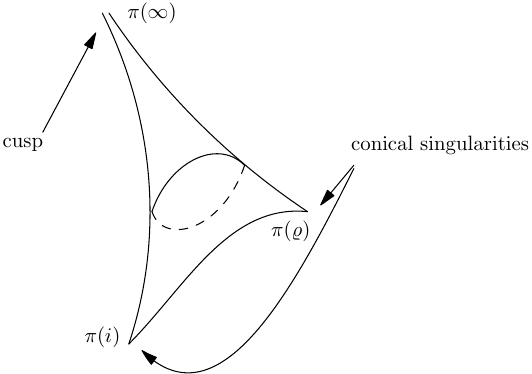} 
\caption{The modular surface~$X=\Gamma\backslash\h$.}\label{fig:modsurface}
\end{figure}

\begin{figure}[h]
\centering
\includegraphics{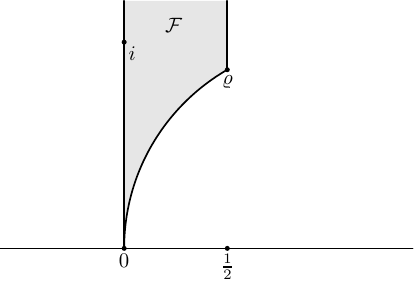} 
\caption{Fundamental domain $\mc F$ for $\Gamma$.}\label{fig:funddom2}
\end{figure}

A model of~$X$ is given by the (closed) \emph{fundamental domain} 
\[
 \mc F_0\ccoloneqq \left\{ z \in \h \ \left\vert\  |z|\geq 1,\ |\Rea z |\leq \tfrac12\right.\right\}
\]
(see Figure~\ref{fig:funddom1}). 
It contains at least one point of any $\Gamma$-orbit, thus $\pi(\mc F_0) = X$. Only points in the boundary of~$\mc F_0$ can be identified under the action of~$\Gamma$, namely the left vertical boundary is mapped to the right one by the element 
\begin{equation}\label{def:T}
 T \ccoloneqq \bmat{1}{1}{0}{1}\,,
\end{equation}
which acts on~$\h$ by $T\act z = z+1$, and the left bottom boundary (the arc from~$\overline{\varrho}$ to~$i$) is mapped to the right bottom boundary (the arc from~$\varrho$ to~$i$) by 
\begin{equation}\label{def:S}
 S \ccoloneqq \bmat{0}{1}{-1}{0}\,,
\end{equation}
which acts on~$\h$ by $S\act z = -1/z$. If we glue  $\mc F_0$ together according to these boundary identifications then we obtain the modular surface~$X$, as illustrated in Figure~\ref{fig:modsurface}. This is just like what we did when we represented $\T = \R/\Z$ as $[0,1]/{\{0=1\}}$. 
Clearly, there is more than one fundamental domain for the modular surface. Another fundamental domain is, e.\,g., 
\[
 \mc F \ccoloneqq \left\{ z\in\h \ \left\vert\ |z-1|\geq 1,\ 0\leq \Rea z\leq \tfrac12    \right.\right\}
\]
(see Figure~\ref{fig:funddom2}). 
It arises from~$\mc F_0$ by cutting off the left half \mbox{$\mc F_L\coloneqq \mc F_0 \cap \{\Rea z < 0\}$} from~$\mc F_0$, gluing $S\act\mc F_L$ to the right half of~$\mc F_0$ and adding all topological boundaries. Thus, 
\[
\mc F \ceqq S\act\overline{\mc F}_L\;\cup\;(\mc F_0\smallsetminus\mc F_L)\,,
\]
where $\overline{\mc F}_L = \mc F_0 \cap \{ \Rea z \leq 0\}$ denotes the closure of~$\mc F_L$ in~$h$. 
For our constructions in Section~\ref{sec:discretization} below, the fundamental set $\mc F$ is more convenient than $\mc F_0$.

The modular surface has an infinite `end' of finite volume, called the \emph{cusp}. In the fundamental domain~$\mc F_0$ it is represented by the strip going to~$\infty$. In terms of~$\Gamma$, the presence of the element~$T$ in~$\Gamma$ caused the presence of this cusp. As we will see, this cusp and the element~$T$ play a special role throughout. 

For completeness we remark that the modular surface is not a hyperbolic surface in the strict sense because it is not a Riemannian manifold but rather an orbifold. It has the two \emph{conical singularities} at~$i$ and~$\varrho$ (see Figures~\ref{fig:funddom1}--\ref{fig:funddom2}). At these points the structure of the quotient space~$X=\Gamma\backslash\h$ is not smooth. The non-smoothness, however, does not influence any step in our discussions.

\subsection{Geometric entity: geodesics.}
\begin{figure}[h]
\centering
\includegraphics{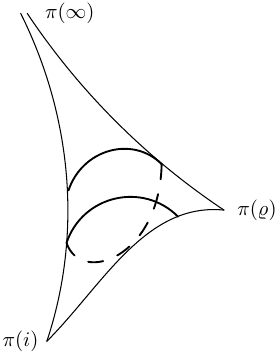} 
\caption{A geodesic on the modular surface.}\label{fig:modgeod}
\end{figure}

Just as in the case of the torus, the `geometric entities' for the modular surface are the periodic geodesics and their lengths.  A geodesic on~$X$ is the image under the projection map~$\pi\colon\h\to X$ of a geodesic on~$\h$, as illustrated in Figure~\ref{fig:modgeod}. 
Geodesics on~$\h$ are infinitely long, but geodesics on~$X$ can be either infinitely long or else periodic and of finite length.  The (primitive) geodesic length spectrum~$L_X$ of~$X$ is by definition the multiset of the lengths of periodic geodesics. The periodic geodesics on~$X$ are closely related to those elements~$g\in\Gamma$ with $|\tr(g)|>2$, the \emph{hyperbolic elements}: For every periodic geodesic~$\wh\gamma$ on~$X$ and any representing geodesic~$\gamma$ of~$\wh\gamma$ on~$\h$ (i.\,e., $\pi(\gamma) = \wh\gamma$) there exists a hyperbolic element~$g\in\Gamma$ such that $g.\gamma$ is a time-shifted version of~$\gamma$, i.\,e., there exists $t_g>0$ such that 
\begin{equation}\label{eq:timeshift}
 g.\gamma(t) = \gamma(t+t_g)\qquad\text{for all $t\in\R$.}
\end{equation}
If in~\eqref{eq:timeshift} the value $t_g$ is minimal among all possible choices of~$g\in\Gamma$, then $g$ is \emph{primitive hyperbolic}. An equivalent characterization is that $g$ is hyperbolic and not of the form~$h^n$ with $h\in\Gamma$ and $n\in\N$. 

Conversely, whenever $\gamma$ is a geodesic on~$\h$ and there exists $g\in\Gamma$ and $t_g>0$ such that~\eqref{eq:timeshift} holds, then $g$ is hyperbolic and $\pi(\gamma)$ is a periodic geodesic on~$X$. Furthermore, every hyperbolic element in~$\Gamma$ time-shifts a unique geodesic on~$\h$. Under this assignment of primitive hyperbolic elements in~$\Gamma$ to periodic geodesics on~$X$, the set of periodic geodesics on~$X$ is bijective to the set of conjugacy classes of the primitive hyperbolic elements in~$\Gamma$, and the (primitive) geodesic length spectrum of~$X$ is the multiset
\[
 L_X \ceqq \left\{ 2\acosh\left(\frac{|\tr(g)|}{2}\right) \ \left\vert\  g\in\HP \vphantom{\Big(\tfrac{|\tr(g)|}{2}\Big)} \right.\right\}_\multiset\,,
\]
where $\HP$ is any set of representatives for the conjugacy classes of primitive hyperbolic elements in~$\Gamma$. The smallest element in $L_X$ is 
\[
 2\acosh\left(\frac{3}{2}\right) \ceqq 2\log\left(\frac{3+\sqrt{5}}2\right)\,,
\]
and more generally the full multiset~$L_X$ consists of \emph{all} numbers of the form  
\[
2 \log\left(\frac{t+\sqrt{t^2-4}}2\right) 
\]
with~$t\in\Z_{\geq 3}$ with multiplicities that can be described explicitly in terms of class numbers of indefinite
quadratic forms.  We refer the interested reader to \cite[Exercises~18-20 in Section~3.7, and the paragraph below them]{Terras} and omit any discussion of this relation here.

\subsection{Spectral entity: Laplace eigenfunctions.} We now introduce the spectral objects we are interested in: the Maass wave forms for $\Gamma$, and the more special Maass cusp forms.

The Laplacian on~$\h$, the \emph{hyperbolic Laplacian}, is 
\begin{equation}\label{eq:hypLaplace}
 \Delta \ccoloneqq -y^2 \big( \partial^2_x + \partial^2_y\big)\qquad (z=x+iy)\,.
\end{equation}
The differential operator~$\Delta$ commutes with all elements of the group~$G=\PSL_2(\R)$ of orientation-preserving Riemannian isometries; the factor~$y^2$ in~\eqref{eq:hypLaplace} corresponds to the factor~$y^{-2}$ in the formula of the line element of the Riemannian metric in~\eqref{eq:lineelement}. Initially, $\Delta$ is defined as an operator on all functions~$\h\to\C$ that are twice partially differentiable. However, it can also be understood as an operator on more general spaces. We refer to~\cite{Hoer1, Venkov} for extensive discussions. 

Now let $u\colon\h\to\C$ be a $\Gamma$-invariant eigenfunction of~$\Delta$, that is, a function satisfying $u(g\act z) = u(z)$ for all~$g\in\Gamma$ and all~$z\in\h$, and
\begin{equation}\label{eq:eigenfunction}
 \Delta u \ceqq s(1-s)u
\end{equation}
for some $s\in\C$. Further below we will see that it is more convenient to work with the \emph{spectral parameter}~$s$ rather than with the eigenvalue~$s(1-s)$ itself. We do not need to specify \emph{a priori} the precise regularity of~$u$, it suffices to require~$u$ to be a hyperfunction or continuous (which is much stronger): since the Laplace operator is elliptic with real-analytic coefficients, the function~$u$ is automatically real-analytic (see~\cite[Theorem~9.5.1]{Hoer1} or~\cite[Theorem~6.33 and its remarks]{Folland}).

The invariance of~$u$ under the element~$T\in\Gamma$ from~\eqref{def:T} shows that $u$ is $1$-periodic, and hence has a Fourier expansion of the form
\[
 u(x+iy) \ceqq \sum_{n\in\Z} a_n(y)\, e^{2\pi i n x}\,.
\]
By separation of variables in~\eqref{eq:eigenfunction} we see that each function~$a_n$ is a solution of a second-order differential equation (depending on~$s$), a modified Bessel differential equation. This equation has two independent solutions, one exponentially big and one exponentially small as~$y\to\infty$, except if~$n=0$, where two independent solutions are~$y^{s}$ and~$y^{1-s}$ for~$s\not=\tfrac12$, and $y^{1/2}$ and $y^{1/2}\log y$ for~$s=\tfrac12$. If we assume in addition that $u$ has polynomial growth at infinity, in which case $u$ is called a \emph{Maass wave form} for~$\Gamma$, then the Fourier expansion becomes
\[
 u(x+iy) \ceqq c_1 y^s + c_2 y^{1-s} + y^{\frac12} \sum_{\substack{n\in\Z\\ n\not=0}} A_n\, K_{s-\frac12}(2\pi |n| y)\, e^{2\pi i n x}\,,
\]
where the first two terms must be replaced by~$c_1 y^{1/2} + c_2 y^{1/2}\log y$ if~$s=\tfrac12$. Here $K_\nu$ is the appropriately normalized solution of the Bessel differential equation that is exponentially small at infinity, the so-called \emph{modified Bessel function of the second kind with index~$\nu\in\C$}, whose precise definition plays no role in our further discussion and is therefore omitted.  The~$A_n$ are complex numbers that automatically have polynomial growth. 

If we further assume that $u$ is bounded, then $c_1 = c_2 = 0$ and 
\[
 u(x+iy) \ceqq y^{\frac12} \sum_{\substack{n\in\Z\\ n\not=0}} A_n\, K_{s-\frac12}(2\pi |n| y)\,e^{2\pi i n x}\,.
\]
In this case, the function~$u$ has rapid decay at infinity and is called a \emph{Maass cusp form with spectral parameter~$s$}. It is known that the real part of~$s$ then always lies between~$0$~and~$1$. Since any Maass wave form~$u$ is $\Gamma$-invariant, we can also consider $u$ as a true function on~$X=\Gamma\backslash\h$, and characterize Maass cusp forms as eigenfunctions of~$\Delta$ on~$X$ having rapid decay as their argument tends to the cusp.

The Friedrichs extension allows us to define~$\Delta$ as an operator on the Hilbert space $L^2(X)$, which can be understood as the space of the (Lebesgue-equivalence classes of) $\Gamma$-invariant functions~$\h\to\C$ that are locally square-integrable~\cite{Venkov}. The \mbox{$L^2$-}eigenfunctions of~$\Delta$ on~$X$ are the constant functions (with eigenvalue~$0$) and the Maass cusp forms, whose eigenvalues are positive and tend to infinity, giving an $L^2$-Laplace spectrum
\[
 \sigma(X) \ceqq \bigl\{ 0,\ 91.141\cdots,\ 148.432\cdots,\ 190.131\cdots,\ \ldots \bigr\}
\]
whose elements can be computed numerically to high precision~\cite{BSV}, but are not known in closed form.

\subsection{Dynamical zeta function.} An analogue of the dynamical zeta function~$\zeta_\T$ of the torus is the Selberg zeta function~$Z_X$, which has an Euler product given by the lengths of periodic geodesics and an Hadamard product in terms of the Laplace resonances (i.e., spectral parameters of generalized eigenfunctions). More precisely, $Z_X(s)$ is defined for~$\Rea s >1$ by 
\[
  Z_X(s) \ceqq \prod_{\ell\in L_X} \prod_{k=0}^\infty \bigl( 1- e^{-(s+k)\ell}\bigl)\,, 
\] 
and the analogue of~\eqref{eq:relation} is Selberg's theorem that this function extends meromorphically to~$\C$ and vanishes if~$s$ is a spectral parameter. See~\cite{Selberg} or~\cite[Chapter~7]{Venkov}.

\section{The cohomological interpretation of Maass cusp forms}\label{sec:mcf}
We now turn to the first step in the passage from geodesics on the modular surface~$X$ to Maass cusp forms for~$\Gamma$: the interpretation of Maass cusp forms in terms of \emph{parabolic $1$-cohomology} as provided in \cite{BLZm}. 

The essential part of this cohomological interpretation, of which we take advantage here, is that every Maass cusp form~$u$ with spectral parameter~$s$ is characterized by a vector~$(c_g^u)_{g\in\Gamma}$ of functions $\mP^1(\R)\to\C$ given by integrals of the form
\begin{equation}\label{def:cg}
 c_g^u(t) \ceqq \int_{g^{-1}\infty}^\infty \omega_s(u,t)
\end{equation}
for~$t\in\R$, and at~$\infty$ by smooth ($C^\infty$) extension (see below for a definition). Here, $\omega_s(u,\cdot)$ is a certain closed \mbox{$1$-}form on~$\h$ defined below and the integration is along any path in~$\h\cup\mP^1(\Q)$ from~$g^{-1}\infty$ to~$\infty$ with at most finitely many points in~$\mP^1(\Q)$ (and which approaches these, say, within a sector). In fact, we usually take a piecewise geodesic path. The functions~$(c_g^u)_{g\in\Gamma}$ satisfy certain relations among each other, so-called cocycle relations, showing that a suitable cohomology theory is the natural home of this setup.

For completeness of exposition and for the convenience of the reader we provide a rather detailed definition of this cohomology (specialized to the modular group~$\Gamma$), even though these details will not be needed further on. Readers who want to proceed faster to the final result are invited to skip the remaining part of this section after having read Theorem~\ref{thm:BLZ}. They should interpret the space~$H^1_\parab(\Gamma;\Vect)$ defined below as a vector space whose elements are equivalence classes of maps from~$\Gamma$ to the space of sufficiently regular functions on~$\mP^1(\R)$, where the notion of `sufficiently regular at~$\infty$' depends on the parameter~$s$. Theorem~\ref{thm:BLZ} then states that the assignment of Maass cusp forms~$u$ with spectral parameter~$s$ to the equivalence classes of the vectors~$(c_g^u)$ is linear and injective, and surjects onto $H^1_\parab(\Gamma;\Vect)$. 

For the detailed description we start with a few preparations. The parabolic cohomology will then be seen to be a refinement of the standard group cohomology in order to account for the cusp of the modular surface and the rapid decay of the Maass cusp forms towards this cusp. The name \emph{parabolic} alludes to the fact that elements in~$G$ that stabilize a single point in~$\mP^1(\R)$, such as~$T$, are called parabolic. 

For any~$s\in\C$, we define an action of~$G$ on partial functions~$\mP^1(\R)\to\C$ by setting
\begin{equation}\label{def:taus}
 \tau_s(g^{-1})f(t) \ccoloneqq \big(g'(t)\big)^s f(g\act t)
\end{equation}
(sometimes also denoted $f\vert_{2s}g$) wherever it is defined. We recall that such a partial function need not be defined on all of~$\mP^1(\R)$. In the situation of~\eqref{def:taus}, the function~$\tau_s(g^{-1})f$ will not be defined on~$g^{-1}\act\infty$ (and maybe additional points).

Let $\Vect$ (called the space of \emph{smooth, semi-analytic vectors of the principal series representation with spectral parameter $s$ in the line model}) denote the space of smooth functions~$\varphi\colon\mP^1(\R)\to\C$ that are real-analytic on~$\R$ up to a finite set that may depend on~$\varphi$, with the action~\eqref{def:taus}. Smoothness at the point~$\infty$ here means that the map
\[
\tau_s(S)\varphi\colon t\mapsto |t|^{-2s} \varphi\big(-\tfrac1t\big) 
\]
extends smoothly ($C^\infty$) to the point~$0$ (recall the element~$S$ from~\eqref{def:S}). For completeness we remark that in~\cite{BLZm} the space~$\Vect$ is denoted~$\mc V_s^{\omega^*,\infty}$.

The vector space $Z^1_\parab(\Gamma;\Vect)$ of \emph{parabolic $1$-cocycles} is then the space of maps~$c\colon\Gamma\to \Vect$ such that 
\begin{itemize}
\item
for all~$g,h\in\Gamma$, we have
 \begin{equation}\label{eq:coc_relation}
  c_{gh} \ceqq \tau_s(h^{-1})c_g + c_h\,,
 \end{equation}
where $c_g$ denotes the function~$c(g)$, and
\item
there exists $\varphi\in \Vect$ such that 
\begin{equation}\label{eq:coc_Phi}
 c_T \ceqq \tau_s(T^{-1})\Phi - \Phi\,.
\end{equation}
(For general discrete subgroups we would need a similar condition for representatives of each conjugacy class of parabolic elements.) 
\end{itemize}
The subspace~$B^1(\Gamma;\Vect)$ of \emph{\mbox{$1$-}coboundaries} consists of the maps $c\colon\Gamma\to \Vect$ for which there exists $\varphi\in\Vect$ such that 
\begin{equation}\label{eq:coboundary}
 c_g \ceqq \tau_s(g^{-1})\varphi -\varphi \qquad\text{for every $g\in\Gamma$.}
\end{equation}
For $c\in B^1(\Gamma;\Vect)$ and $\varphi\in\Vect$ as in~\eqref{eq:coboundary} we find for all $g,h\in\Gamma$ the identity
\begin{align*}
 c_{gh} & = \tau_s\big( (gh)^{-1} \big)\varphi - \varphi = \tau_s\big( h^{-1}g^{-1} \big)\varphi - \varphi 
 \\
 & = \tau_s\big(h^{-1}\big) \Big( \tau_s\big(g^{-1}\big)\varphi  - \varphi\Big) + \tau_s\big(h^{-1}\big)\varphi - \varphi
 \\
 & = \tau_s\big(h^{-1}\big)c_g + c_h\,,
\end{align*}
which shows that every $1$-coboundary is a $1$-cocycle, and also that  $B^1(\Gamma;\Vect)$ is a subspace of~$Z^1_\parab(\Gamma;\Vect)$ (take $\Phi=\varphi$ in~\eqref{eq:coc_Phi}). The quotient space
\[
 H^1_\parab(\Gamma;\Vect) \ccoloneqq Z^1_\parab(\Gamma;\Vect)/B^1(\Gamma;\Vect)
\]
is called the \emph{space of parabolic $1$-cohomology classes} with values in~$\Vect$. 

For any two real-analytic functions~$u,v$ on~$\h$ we define the \emph{Green's form} to be the real-analytic \mbox{$1$-}form
\[
 [u,v] \ccoloneqq \frac{\partial u}{\partial z} \cdot v\cdot dz \cplus u\cdot\frac{\partial v}{\partial {\overline z}}\cdot d\overline z\,,
\]
which is easily seen to be closed (i.\,e., $d[u,v]=0$) if $u$ and~$v$ are eigenfunctions of~$\Delta$ with the same eigenvalue. For any~$s\in\C$ and any~$t\in\R$ the function~$R(t;\cdot)^s\colon \h\to\C$, where
\[
 R(t;z)\ccoloneqq \Ima\frac{1}{t-z}\,,
\]
is an eigenfunction of~$\Delta$ with eigenvalue~$s(1-s)$. Therefore, if $u$ is a Maass cusp form with spectral parameter~$s$, then for any~$t\in\R$ the $1$-form 
\[
 \omega_s(u,t) \ccoloneqq \big[u, R(t;\cdot)^s\big]
\]
is closed. From this it follows that, for any~$g\in\Gamma$, the integral 
\begin{equation}\label{eq:paramintegral}
c_g^u(t) \ccoloneqq \int_{g^{-1}\infty}^\infty \omega_s(u,t) 
\end{equation}
is independent of the chosen path from~$g^{-1}\infty$ to~$\infty$. The integral is convergent due to the rapid decay of~$u$ at the cusp. The regularities of~$u$ and~$R(\cdot\,;\cdot)^s$ yield \mbox{$c^u_g\in \Vect$}. Furthermore, the \mbox{$\Gamma$-}invariance of~$u$ implies the transformation formula 
\begin{equation}\label{eq:transform}
 \tau_s(g)\int_a^b \omega_s(u,t) \ceqq \int_{g\cdot a}^{g\cdot b} \omega_s(u,t) \qquad (g\in\Gamma,\ a,b\in\mP^1(\R),\ t\in\R)
\end{equation}
and from this one easily deduces that the map~$c^u$ satisfies the cocycle relation~\eqref{eq:coc_relation} and the relation in~\eqref{eq:coc_Phi} and hence is a parabolic cocycle. Then we have:

\begin{thm}[\cite{LZ01,BLZm}]\label{thm:BLZ}
For $s\in\C$, $\Rea s\in (0,1)$, the map~$u\mapsto [c^u]$ defines an isomorphism of vector spaces 
\[
 \{\text{Maass cusp forms with  spectral parameter $s$}\}\, \overset{\thicksim}{\longrightarrow}\, H^1_\parab(\Gamma;\Vect)\,.
\]
\end{thm}

\section{Discretization of geodesics}\label{sec:discretization}

In this section we will discuss the second step in the passage from geodesics on the modular surface~$X$ to Maass cusp forms for~$\Gamma$: the construction of  a discretization of the motion along the geodesics on~$X$. 

The two elements (generators)
\begin{equation}\label{def:TT}
 T_1 \coloneqq \mat{1}{0}{1}{1}\qquad\text{and}\qquad T_2 \coloneqq \mat{1}{1}{0}{1}
\end{equation}
of~$\Gamma$  and the map 
\[
 F\colon (0,\infty)\smallsetminus\Q \to (0,\infty)\smallsetminus\Q
\]
given by the two branches
\begin{equation}\label{def:F}
\begin{cases}
 (0,1)\smallsetminus\Q \stackrel{\sim}{\longrightarrow} (0,\infty)\smallsetminus\Q\,, & x\mapsto T_1^{-1} x = \frac{x}{1-x}
 \\[2mm]
 (1,\infty)\smallsetminus\Q \stackrel{\sim}{\longrightarrow} (0,\infty)\smallsetminus\Q\,, & x\mapsto T_2^{-1} x = x-1
\end{cases} 
\end{equation}
will play a crucial role. By iterating the map~$F$ we get a \emph{discrete(-time) dynamical system} 
\begin{equation}\label{def:Fdyn}
 \N_0 \times \bigl((0,\infty)\smallsetminus\Q\bigr) \to (0,\infty)\smallsetminus\Q\,,\quad (n,x) \mapsto F^n(x)\,,
\end{equation}
which we denote for short by~$F$ as well. (It will always be clear if~$F$ refers to the map in~\eqref{def:F} or to the map in~\eqref{def:Fdyn}.) We will show that this discrete dynamical system can be thought of as a discrete version of the geodesic flow on~$X$: The map~$F$ and its iterates capture the essential geometric and dynamical properties of the geodesic flow that will be needed for establishing the relation between the geodesics on $X$ and the Maass cusp forms for~$\Gamma$. In particular, the orbits of the map~$F$ describe the future behavior of (almost all) geodesics on~$X$, and periodic geodesics on~$X$ correspond to points~$x\in(0,\infty)\smallsetminus\Q$ with periodic (i.\,e., finite) orbits under~$F$.\footnote{We remark that the formula for~$F$ is identical to the map~$\Phi$ given in \cite[Section~1.1, Lemma]{Choie_Zagier} in connection with the so-called rational period functions.} 

The construction of~$F$ from the geodesic flow on~$X$ proceeds in several steps: We first choose a `good' cross section (in the sense of Poincar\'e) for the geodesic flow on~$X$, i.\,e., a subset~$\wh C$ of the unit tangent bundle of~$X$ that is intersected by all periodic geodesics at least once, and each intersection between any geodesic on~$X$ and~$\wh C$ is discrete. We refer to the discussion below for precise definitions. The choice of~$\wh C$ yields a first return map, which is the map that assigns to each element $\wh v\in\wh C$ the next intersection between~$\wh C$ and the geodesic on~$X$ starting at time~$0$ in the direction~$\wh v$. The first return map provides a first discretization of the geodesic flow on~$X$. 

Then we choose a `good' set of representatives for~$\wh C$, i.\,e., a subset~$C^*$ of the unit tangent bundle of $\h$ that is bijective to~$\wh C$ with respect to the canonical quotient map. The specific properties of~$C^*$ will allow us to semi-conjugate the first return map to a map on $(0,\infty)\smallsetminus\Q$, which is precisely the map~$F$. 

The construction we will present below is a special case of the algorithm in~\cite{Pohl_Symdyn2d} for finding good discretizations for geodesic flows on much more general hyperbolic surfaces. We refer to~\cite{Pohl_Symdyn2d} for further details and all omitted proofs, in particular to~\cite[Proposition~8.2, Theorem~8.15, Corollary~8.16]{Pohl_Symdyn2d} and their specialization to the modular surface as in~\cite[Example~3.3]{Pohl_Symdyn2d}.

As in Section~\ref{sec:discretization}, readers who want to proceed faster to the final result are invited to skip the remaining part of this section after having read Theorem~\ref{thm:CS} below. In Section~\ref{sec:TO} only the map~$F$ will be needed, not the details of its construction.

\subsection{Geodesics.} While in Section~\ref{sec:mod_surface} we used the notion of geodesics in the sense~\ref{geod2} (adapted to the hyperbolic plane and the modular surface in place of the real line and the torus), we now also need geodesics in the sense~\ref{geod1}. 

A geodesic~$\gamma$ on~$\h$ in the sense~\ref{geod1} is completely determined by requiring that it passes through a given point~$z\in\h$ at time~$t=0$ in a given direction. Recall that we consider only geodesics of unit speed, so that the speed in the given direction does not form another parameter. Therefore we may identify geodesics in the sense~\ref{geod1} with the set of all unit length direction vectors at all points of~$\h$, thus, with the \emph{unit tangent bundle}~$S\h$ of~$\h$. 

For $v\in S\h$ we let $\gamma_v\colon\R\to\h$ be the (unique) geodesic on~$\h$ such that 
\begin{equation}\label{def:geodH}
 \gamma_v'(0)=v\,.
\end{equation}
Both the tangent vector~$\gamma_v'(0)$ to~$\gamma_v$ at time~$t=0$ and  the element~$v\in S\h$ are combinations of position and direction, the position~$\gamma_v(0)$ being the \emph{base point}~$\base(v)\in\h$. The \emph{geodesic flow} on~$\h$ (the motion along geodesics on~$\h$) is the map
\begin{equation}\label{def:geodflowH}
\R\times S\h \to S\h\,,\quad (t,v)\mapsto \gamma_v'(t)\,.
\end{equation}
The action of~$G$ on~$\h$  by Riemannian isometries induces an action of~$G$ on~$S\h$  by 
\[
 g\act v \ccoloneqq (g\act \gamma_v)'(0) \qquad (g\in G,\ v\in S\h)\,.
\]
The \emph{unit tangent bundle} of~$X$ is then just the quotient
\[
 SX \ceqq \Gamma\backslash S\h\,.
\]
We denote the projection map
\begin{equation}\label{def:projmaptangent}
\pi\colon S\h\to SX
\end{equation}
with the same symbol as the projection map~$\h\to X$ from~\eqref{def:projmap}. The context always clarifies which one is meant. We typically denote a geodesic on $\h$ by $\gamma$ and a unit tangent vector in $S\h$ by $v$, and use $\wh\gamma$ and $\wh v$ for the corresponding geodesic $\pi(\gamma)$ on $X$ and unit tangent vector $\pi(v)\in SX$. In analogy with \eqref{def:geodH}, for any $\wh v\in SX$ we let $\wh\gamma_v$ denote the geodesic on $X$ determined by 
\[
\wh\gamma_v'(0) = \wh v\,.
\]
Also the \emph{geodesic flow} on~$X$ is inherited from the geodesic flow on~$\h$ as defined in~\eqref{def:geodflowH}, and hence is the map
\[
 \R\times SX\to SX\,,\quad (t,\wh v)\mapsto \wh\gamma_v'(t)\,.
\]

\subsection{Cross section.} 
By a \emph{cross section} we mean (slightly deviating from the standard definition) a subset~$\wh C$ of~$SX$ such that 
\vspace*{-.5mm}
\begin{enumerate}[label=$\mathrm{(\mathbf{C\arabic*})}$, ref=$\mathrm{(\mathbf{C\arabic*})}$]
 \item every periodic geodesic on~$X$ intersects $\wh C$. In other words, for any periodic geodesic~$\wh\gamma$ there exists~$t\in\R$ such that $\wh\gamma'(t) \in \wh C$.
 \item each intersection of any geodesic on~$X$ with~$\wh C$ is discrete. In other words, for any geodesic~$\wh\gamma$ and $t\in\R$ with $\wh\gamma'(t)\in\wh C$ there exists~$\eps>0$ such that 
 \[
  \wh\gamma'\big((t-\eps, t+\eps)\big) \cap \wh C = \big\{\wh\gamma'(t)\big\}\,.
 \]
\end{enumerate}
We  define a \emph{set of representatives}~$C^*$ for a cross section~$\wh C$ to be a subset of~$S\h$  that is bijective to~$\wh C$ under the projection map~$\pi$ from \eqref{def:projmaptangent}. (We write~$C^*$ rather than~$C$ because the latter traditionally denotes the full preimage of~$\wh C$ in~$S\h$.) Of course, to characterize a cross section~$\wh C$ it suffices to provide a set of representatives, but choosing a cross section and a set of representatives that serves our purposes is an art. For the modular surface we will take 
\[
 C^* \ccoloneqq \{ v\in S\h \mid \base(v)\in i\R^+,\ \gamma_v(\infty) \in (0,\infty)\smallsetminus\Q \}
\]
as set of representatives, where 
\[
 \gamma_v(\infty) \ccoloneqq \lim_{t\to\infty} \gamma_v(t)\,.
\]
\begin{figure}[h]
\centering
\includegraphics{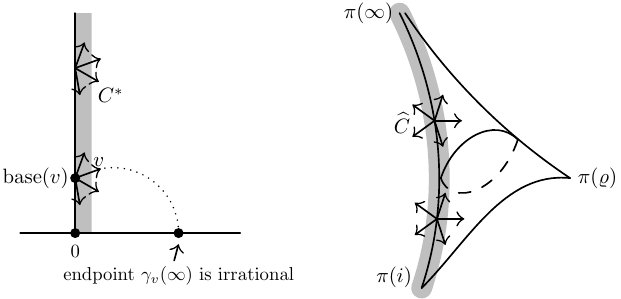} 
\caption{The set of representatives~$C^*$ and the cross section~$\wh C$. The gray shadows indicate the directions of the elements of~$\wh C$ and~$C^*$.}\label{fig:crosssection}
\end{figure}

\noindent
The associated cross section 
\[
 \wh C \ccoloneqq \pi(C^*)
\]
is the set of unit tangent vectors~$\wh v\in SX$ sitting on the geodesic from~$\pi(i)$ to~$\pi(\infty)$ such that the geodesic emanating from~$\wh v$ does not converge to the cusp~$\pi(\infty)$ in future or past time. A pictorial representation of~$C^*$ and~$\wh C$ is given in Figure~\ref{fig:crosssection}. Choosing a set of representatives~$C^*$ such that the base points of its elements forms the geodesic from~$0$ to~$\infty$ in~$h$ is motivated by the integral expression in~\eqref{eq:paramintegral}. Its effect will become clearer in Section~\ref{sec:TO}.

\subsection{Discretization.} 
We will now show how to relate the geodesic flow on~$X$ to a discrete dynamical system on (a subset of)~$\R_{>0}$. In the case of the modular surface, this construction is closely related to continued fractions, more precisely to Farey fractions. The reader interested in this connection may find the articles \cite{Artin,Richards,Series,Katok_Ugarcovici} useful. 

\begin{figure}[h]
\centering
\includegraphics{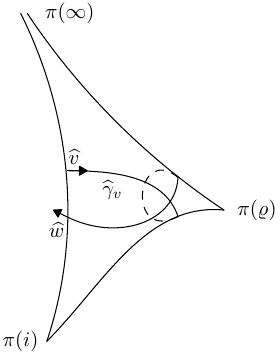} 
\caption{The geodesic determined by $\wh v$ and its first return to $\wh C$.}\label{fig:cross2}
\end{figure}

\begin{figure}[h]
\centering
\includegraphics{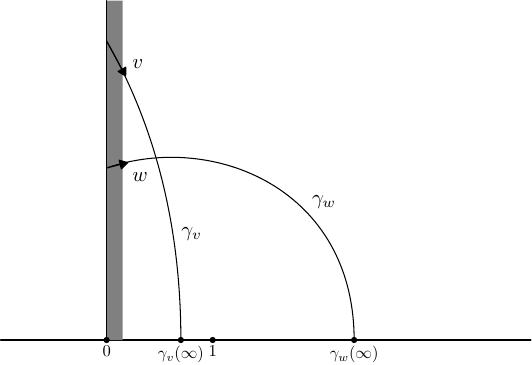} 
\caption{Associated geodesics on $\h$.}\label{fig:next}
\end{figure}

Let $\wh v\in \wh C$ be an element of the cross section and consider the associated geodesic~$\wh\gamma_v$ on~$X$. By the choice of~$\wh C$, the geodesic~$\wh\gamma_v$ intersects~$\wh C$ again in future time. Let $t_0>0$, the \emph{first return time}, be the minimal positive number such that 
\[
 \wh w\ccoloneqq\wh\gamma_v'(t_0) \in \wh C.
\]
(See Figure~\ref{fig:cross2}.) 
Let $v,w\in C^*$ be the elements in the set of representatives corresponding to $\wh v, \wh w$, and $\gamma_v, \gamma_w$ the associated geodesics on $\h$.  
(See Figure~\ref{fig:next}.) 
Since the unit tangent vector $\gamma_v'(t_0) \in S\h$ projects to $\wh w$ under $\pi$, that is,
\[
 \pi\big( \gamma_v'(t_0) \big) \ceqq \wh w\,,
\]
there exists a unique element $g\in\Gamma$ such that 
\[
 \gamma_v'(t_0) \ceqq g\act w\,.
\]
This element is characterized by the property that
\begin{equation}\label{eq:nextelem}
 \gamma_v'(t_0) \;\in\; g\act C^*,
\end{equation}
i.\,e., by the first intersection of~$\gamma_v$ with some $\Gamma$-translate of~$C^*$ after passing through~$v=\gamma_v'(0)$.  To find the element~$g$ we consider the neighboring translates of the fundamental domain~$\mc F$ and the relevant translates of~$C^*$.
\begin{figure}[h]
\centering
\includegraphics{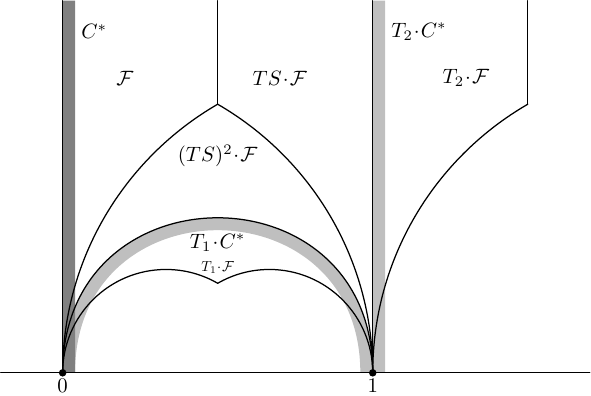} 
\caption{Relevant $\Gamma$-translates of~$\mc F$ and~$C^*$.}\label{fig:forward}
\end{figure}

We observe that, as shown in Figure~\ref{fig:forward}, the unit tangent vector $\gamma_v'(t_0)$ can be only in $T_1\act C^*$ or $T_2\act C^*$ with $T_1,T_2$ as in~\eqref{def:TT}. 
\begin{figure}[h]
\centering
\includegraphics{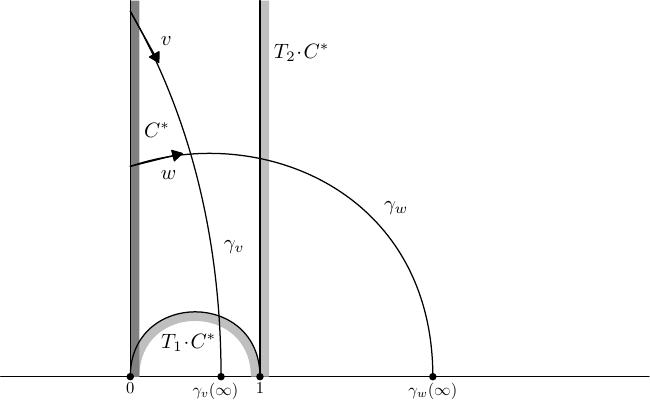} 
\caption{Next intersection. Recall that $v=\gamma'_v(0)$ and $w=\gamma'_w(0)$.}\label{fig:nextex}
\end{figure}
In Figure~\ref{fig:nextex} we have $g=T_1$, so that here 
\[
w \ceqq T_1^{-1}\gamma_v'(t_0)\,, \qquad  \gamma_w(\infty) \ceqq T_1^{-1}\gamma_v(\infty)\,.
\]
We further observe  that for every point $x\in (0,\infty)\smallsetminus\Q$, no matter which \mbox{$v\in C^*$} with $\gamma_v(\infty)=x$ we consider, we find the same value for the element $g\in\Gamma$ in~\eqref{eq:nextelem}. This is caused by the property of~$C^*$ that for~$j\in\{1,2\}$ the set of base points of the vectors in~$T_j\act C^*$ split the hyperbolic plane~$\h$ into two half-spaces and that $T_j\act C^*$ consists of \emph{all relevant} vectors pointing into one of these half-spaces. Therefore the element~$g$ defined by~\eqref{eq:nextelem} depends only on~$x$, not on the specific element~$v\in C^*$ with~$\gamma_v(\infty)=x$. The procedure just described induces a \emph{discrete dynamical system} 
\begin{equation}\label{def:F2}
 F\colon (0,\infty)\smallsetminus\Q \to (0,\infty)\smallsetminus\Q\,,
\end{equation}
where for each $x\in (0,\infty)\smallsetminus\Q$, we pick $v\in C^*$ such that $\gamma_v(\infty)=x$, let $g$ be the element in $\Gamma$ such that $\gamma_v'(t_0)\in g\act C^*$ and set
\[
 F(x) \ccoloneqq g^{-1}\act x\,.
\]

\begin{thm}[\cite{Pohl_Symdyn2d}]\label{thm:CS}
The set $\wh C$ is a cross section for the geodesic flow on $X$, and $C^*$ is a set of representatives for $\wh C$. The induced discrete dynamical system (as in~\eqref{def:F2}) is the map~$F$ as given in~\eqref{def:F}.
\end{thm}

\section{Transfer operators and Maass cusp forms}\label{sec:TO}

In this section we carry out the third and final step in the passage from geodesics on the modular surface~$X$ to Maass cusp forms for~$\Gamma$: to tie together the discrete dynamical system~$F$ from Section~\ref{sec:discretization} and the cohomological interpretation of Maass cusp forms from Section~\ref{sec:mcf}.

The mediating object between both sides is the \emph{transfer operator family}~$(\TO_s)_{s\in\C}$ associated to~$F$.  
The \emph{transfer operator}~$\TO_s$ with parameter~$s$ acts on the vector space of functions from~$(0,\infty)$ to~$\C$ and is given by 
\begin{equation}\label{def:TO}
 \TO_sf(t) \ccoloneqq \sum_{w\in F^{-1}(t)} |F'(w)|^{-s} f(w)
\end{equation}
for $f\in\C^{(0,\infty)}$, $t\in (0,\infty)$. 
This operator has its origin in the thermodynamic formalism of statistical mechanics. It is a generalization of the transfer matrix for \mbox{lattice--spin} systems, which is used to find equilibrium distributions. The weight, being the $(-s)$th-power of the derivative of~$F$, is motivated within this framework, where~$s$ serves as an inverse Boltzmann constant and temperature. From a purely mathematical point of view, this operator can be seen as an evolution operator or as a graph Laplacian on a somewhat generalized graph, in both cases with appropriate weights. 
The explicit expression for~$F$ allows us to evaluate~\eqref{def:TO} in our special case to
\[
 \TO_sf(t) \ceqq f(t+1) \cplus (t+1)^{-2s} f\Bigl(\frac{t}{t+1}\Bigr)\,,\qquad t>0\,,
\]
or, using \eqref{def:taus}, to
\[
 \TO_s \ceqq \tau_s(T_1^{-1}) \cplus \tau_s(T_2^{-1})\,.
\]
(This simple formula is for the modular group only. For other groups one can have a vector of more complicated finite sums.)

The correspondence that we have been aiming at is a bijection between the eigenfunctions of~$\TO_s$ with eigenvalue~$1$ and the Maass cusp forms with spectral parameter~$s$. More precisely, we have the following theorem.

\begin{thm}[\cite{Moeller_Pohl, Pohl_mcf_general}]
Let $s\in\C$, $1>\Rea s>0$. Then for any Maass cusp form~$u$ with spectral parameter~$s$, the function~$f_u\colon (0,\infty)\to\C$ defined by
\begin{equation}\label{eq:f_path}
 f_u(t) \ccoloneqq \int_0^\infty \omega_s(u,t)
\end{equation}
is a real-analytic eigenfunction of~$\TO_s$ with eigenvalue~$1$. The map $u\mapsto f_u$ is a linear isomorphism between the space of Maass cusp forms with spectral parameter~$s$ and the space of real-analytic eigenfunctions~$f$ of~$\TO_s$ with eigenvalue~$1$ for which the map $\R\smallsetminus\{0\} \to \C$ defined by 
\begin{equation}\label{eq:mapdecay}
 \begin{cases}
  f & \text{on $(0,\infty)$}
  \\
  -\tau_s(S)f & \text{on $(-\infty,0)$}
 \end{cases}
\end{equation}
extends smoothly to~$0$.
\end{thm}

\begin{figure}[h]
\centering
\includegraphics{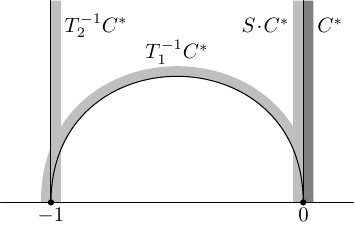} 
\caption{Relevant $\Gamma$-translates for proof of Theorem.}\label{fig:backwards}
\end{figure}

We will now explain the main steps of the proof with an emphasis on intuition and heuristics. Some steps will be omitted, most prominently some discussions of convergence and regularities. We hope to convince the reader that a major part of the proof is encoded in Figure~\ref{fig:backwards} and that the choice of the integral path in \eqref{eq:f_path} and the function in \eqref{eq:mapdecay} is natural. 

\noindent
\textbf{Proof (key elements).} We present the main ideas of the proof, split into four steps.

\textbf{Step 1: Relation between~$\TO_s$ and~$C^*$.}
We first reconsider the transfer operator~$\TO_s$ and its domain~$\C^{(0,\infty)}$. 
We may think of any~$f\in\C^{(0,\infty)}$ as being a mass distribution or density on $(0,\infty)$ of which the transfer operator evaluates its \mbox{$s$-}weighted evolution under one application of~$F$. Recalling that $F$ is a discrete version of the geodesic flow on~$X$, that $\TO_s$ is a weighted evolution operator of~$F$, and that the essential ingredient of this discretization is the set~$C^*$, we may intuitively think of~$f$ as being a `shadow' of some function~$f^*$ on~$C^*$ that is constant on any set of the form
\[
 E_t\ccoloneqq \{ v\in C^* \mid \gamma_v(\infty) = t\}\qquad (t\in (0,\infty))\,.
\]
Thus, 
\[
 f(t) \ceqq f^*(v) \qquad\text{for any~$v\in E_t$.}
\]
When developing the formula for~$F$ we asked where the geodesics determined by the elements in~$C^*$ go to. In the expression for~$\TO_s$, the preimage of~$F$ is used. Hence, when building $\TO_s$, we may alternatively ask where these geodesics come from. For the modular group~$\Gamma$, the relevant sets are~$T_1^{-1}C^*$ and~$T_2^{-1}C^*$. (See Figure~\ref{fig:backwards}.) 

\textbf{Step 2: Relation between Maass cusp forms and~$C^*$.} Let $u$ be a Maass cusp form with spectral parameter~$s$. We use the characterization of~$u$ via a cocycle class in the space~$H^1_\parab(\Gamma;\Vect)$ from Theorem~\ref{thm:BLZ}, and then use the family of functions~$(c_g^u)_{g\in\Gamma}$ from \eqref{eq:paramintegral} as a representative for this cocycle class. We think of each~$c_g$ as being the integral along the geodesic from~$g^{-1}\infty$ to~$\infty$, or even better, as an integral over the set of unit tangent vectors to this geodesics. In particular, for~$g=S$ we have $S^{-1}\infty=0$, so that
\begin{equation}\label{eq:cS}
 c_S^u(t) \ceqq \int_0^\infty \omega_s(u,t)\qquad (t\in\R)
\end{equation}
is the integral along the geodesic from~$0$ to~$\infty$. Thus, in an intuitive way, we may think of~$c_S^u$ as an integral over~$C^*\cup S\act C^*$ and of each value~$c_S^u(t)$ as the mean of some (fictive) function~$\tilde c_S^u(t)$ defined on~$C^*\cup S\act C^*$. 

\textbf{Step 3: From Maass cusp forms to eigenfunctions of $\TO_s$.} Let $u$ be a Maass cusp form with spectral parameter $s$ and $(c_g^u)_{g\in\Gamma}$ the associated family of functions from~\eqref{eq:paramintegral}. We want to associate to~$u$ in a natural way an eigenfunction~$f$ of~$\TO_s$ with eigenvalue~$1$. The intuitive way of thinking of~$c^u_S$ and any function~$f$ as objects related to~$C^*$ suggests using~$C^*$ as mediating element. Staying with this intuition, we should restrict $c_S^u$ to an integral over~$C^*$ and use a relation like $f^*(v) = \tilde c_S^u(t)\vert_{C^*}$ for~$v\in E_t$. In terms of the actual objects (and their rigorous definitions) we are led to set
\begin{equation}\label{eq:f_heur}
 f \ccoloneqq c_S^u\vert_{(0,\infty)}\,,
\end{equation}
which is precisely~\eqref{eq:f_path}. 

We now show that \eqref{eq:f_heur} indeed defines an eigenfunction of $\TO_s$ with eigenvalue~$1$. So far we have used in~\eqref{eq:cS}, and hence in~\eqref{eq:f_heur},  the geodesic from~$0$ to~$\infty$ as path of integration. Since the \mbox{$1$-}form~$\omega_s(u,t)$ is closed, we may change the path to be the geodesic from~$0$ to~$-1$ followed by the geodesic from~$-1$ to~$\infty$:
\[
 \int_0^\infty \omega_s(u,t) \ceqq \int_0^{-1} \omega_s(u,t) \cplus \int_{-1}^\infty \omega_s(u,t)\,.
\]
Using the transformation formula~\eqref{eq:transform} we now find, for any~$t\in (0,\infty)$, 
\begin{align*}
 f(t) & \ceqq  \int_0^\infty \omega_s(u,t) 
 \\
 & \ceqq \int_{T_1^{-1}0}^{T_1^{-1}\infty} \omega_s(u,t) \cplus \int_{T_2^{-1}0}^{T_2^{-1}\infty}\omega_s(u,t)
 \\
 & \ceqq \tau_s(T_1^{-1}) \int_0^\infty \omega_s(u,t) \cplus  \tau_s(T_2^{-1}) \int_0^\infty \omega_s(u,t)
 \\
 & \ceqq \tau_s(T_1^{-1})f(t) \cplus \tau_s(T_2^{-1})f(t)\,.
\end{align*}
Therefore $f=\TO_sf$.

\textbf{Step 4: From eigenfunctions of~$\TO_s$ to Maass cusp forms.} Conversely, let $f$ be a real-analytic eigenfunction of~$\TO_s$ with eigenvalue~$1$ that satisfies the requirement in~\eqref{eq:mapdecay}. We want to associate to~$f$ a Maass cusp form~$u$ in a way which inverts the mapping from Step~3 and which is also natural. 
Instead of trying to do this directly, we will define a parabolic $1$-cocycle~$c=c^f$ in $Z^1_\parab(\Gamma;\Vect)$. Theorem~\ref{thm:BLZ} then implies that the cocycle~$c$ is indeed of the form~$c=c^u$ for a unique Maass cusp form~$u$. 

In order to define~$c$ we prescribe it on the group elements~$T$ and~$S$. Applying~\eqref{eq:paramintegral} for $g=T$, in which case the integral in~\eqref{eq:paramintegral} vanishes, motivates setting 
\[
 c_T \ccoloneqq 0\,.
\]
Further, the intuition explained above suggests defining 
\begin{equation}\label{eq:cs_from_f}
c_S \ccoloneqq 
\begin{cases}
f & \text{on $(0,\infty)$}
\\
-\tau_s(S)f & \text{on $(-\infty,0)$\,.}
\end{cases}
\end{equation}
The minus sign in the second row is motivated by the fact that~$S$ `changes the direction' of the geodesic from~$0$ to~$\infty$.  Since the functions in~\eqref{eq:cs_from_f} and~\eqref{eq:mapdecay} coincide, the regularity properties of~$f$ imply that $c_S$ as defined in~\eqref{eq:cs_from_f} on $\mP^1\R\smallsetminus\{0,\infty\}$ extends smoothly to~$0$ and~$\infty$. 

Since $T$ and~$S$ generate all of~$\Gamma$, the cocycle relation~\eqref{eq:coc_relation} dictates the value of~$c$ on all other elements. It remains to show that $c$ is well-defined, which here means that if a word in~$T$, $T^{-1}$ and~$S$ equals the identity in~$\Gamma$, then the corresponding $\Z[\Gamma]$-combination of~$c_T$ and~$c_S$ vanishes. To that end we use the presentation
\[
 \Gamma \ceqq \left\langle\ S,\;T \ \left\vert\ S^2 = \big(T^{-1}S\big)^3 = \id\ \right.\right\rangle 
\]
and show that 
\[
 \tau_s(S)c_S\cplus c_S \quad\text{and}\quad \big(\tau_s\big((ST)^2\big) + \tau_s(ST) + 1\big)\big(\tau_s(S)c_{T^{-1}} + c_{S}\big)
\]
vanish identically. For the first expression, this follows immediately from~\eqref{eq:cs_from_f}. For the second expression we use~$c_T=0$, deduce first $c_{T^{-1}} = -\tau_s(T)c_T=0$ and then find 
\begin{align*}
 \big(\tau_s\big((ST)^2\big) & + \tau_s(ST) + 1\big)\big(\tau_s(S)c_{T^{-1}} + c_{S}\big)
 \\
 & \ceqq \tau_s\big((ST)^2\big)c_S \cplus \tau_s(ST)c_S \cplus c_S 
 \\
 & \ceqq 
 \begin{cases}
  -\tau_s\big(T_2^{-1}\big)f - \tau_s\big(T_1^{-1}\big)f + f & \text{on $(0,\infty)$}
  \\[2mm]
  \tau_s\big(T_1^{-1}S\big) \left[ -\tau_s\big(T_1^{-1}\big)f + f - \tau_s\big(T_2^{-1}\big)f\right] & \text{on $(-1,0)$}
  \\[2mm]
  \tau_s\big(T^{-1}S\big) \left[ f - \tau_s\big(T_2^{-1}\big)f - \tau_s\big(T_1^{-1}\big)f\right]  & \text{on $(-\infty,-1)$\,,}
 \end{cases}
\end{align*}
which vanishes since $f=\TO_sf$. This calculation 
\begin{figure}[H]
\centering
\includegraphics{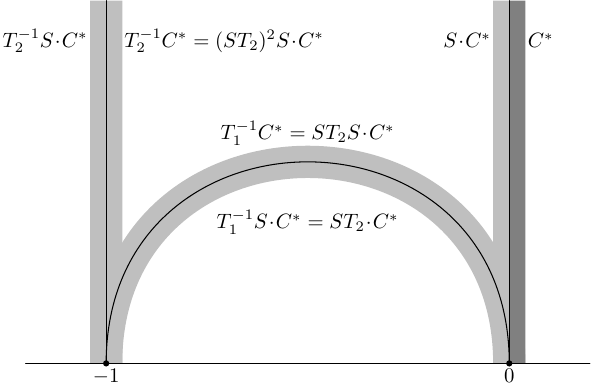} 
\caption{Relevant $\Gamma$-translates for proof of Theorem.}\label{fig:readoff}
\end{figure}
\noindent
can also be read off from Figure~\ref{fig:readoff}, as the reader can verify.
\qed

\section{Recapitulation and closing comments}\label{sec:recap}

We have surveyed an intriguing relation between the periodic geodesics on the modular surface $X=\Gamma\backslash\h$ (`classical mechanical objects') and the Maass cusp forms for $\Gamma$ (`quantum mechanical objects'). For this, we started simultaneously on both ends:

On the geometric side, we developed a discrete version of the (periodic part of the) geodesic flow on the modular surface by means of a cross section in the sense of Poincar\'e. We realized this discretization as a discrete dynamical system on~$(0,\infty)$ by using a well-chosen representation of the cross section on the upper half plane. This step turns the geodesic flow into a discrete and somehow finite object while preserving its essential dynamical features. 

On the spectral side, we characterized the Maass cusp forms as cocycle classes in a certain precise cohomology space. The isomorphism from Maass cusp forms to cocycle classes is given by an integral transform, where a certain \mbox{$1$-}form is integrated along certain geodesics. Even though the cocycle classes remain objects of quantum mechanical nature, this characterization of Maass cusp forms constitutes a first and very important step towards the geometry and dynamics of the modular surface.

Connecting these two sides is the family of transfer operators, which from their definition are purely classical mechanical objects but which clearly exhibit a quantum mechanical nature. These transfer operators depend heavily on the choice of the discretization. The proof of the isomorphism between eigenfunctions of the transfer operators and the parabolic \mbox{$1$-}cocycles clearly shows that the shape of the set of representatives is crucial. Here, it is the set of (almost) all unit tangent vectors that are based on the geodesic from~$0$ to~$\infty$ and that point `to the right'. 

This set of representatives and its $\Gamma$-translates can be seen as a geometric realization of the cohomology. The transfer operator then encodes the cocycle relation. An eigenfunction with eigenvalue~$1$ of the transfer operator obeys a geometric variant of the cocycle relation, and hence can be related to an actual cocycle, which in turn characterizes a Maass cusp form.

\bibliography{PZbib}

\providecommand{\bysame}{\leavevmode\hbox to3em{\hrulefill}\thinspace}
\providecommand{\MR}{\relax\ifhmode\unskip\space\fi MR }
\providecommand{\MRhref}[2]{%
  \href{http://www.ams.org/mathscinet-getitem?mr=#1}{#2}
}
\providecommand{\href}[2]{#2}
\begin{thebibliography}{10}

\bibitem{Artin}
E.~Artin, \emph{Ein mechanisches {S}ystem mit quasiergodischen {B}ahnen}, Abh.
  Math. Sem. Univ. Hamburg \textbf{3} (1924), 170--175.

\bibitem{Bergeron}
N.~{Bergeron}, \emph{{The spectrum of hyperbolic surfaces}}, Les Ulis: EDP
  Sciences; Cham: Springer, 2016.

\bibitem{BSV}
A.~{Booker}, A.~{Str\"ombergsson}, and A.~{Venkatesh}, \emph{{Effective
  computation of Maass cusp forms}}, {Int. Math. Res. Not.} \textbf{2006}
  (2006), no.~12, 34, Id/No 71281.

\bibitem{Boothby}
W.~{Boothby}, \emph{{An introduction to differentiable manifolds and Riemannian
  geometry. 2nd ed}}, {Pure and Applied Mathematics, 120. Academic Press,
  Inc.}, 1986.

\bibitem{Bruggeman_lewiseq}
R.~{Bruggeman}, \emph{{Automorphic forms, hyperfunction cohomology, and period
  functions}}, {J. Reine Angew. Math.} \textbf{492} (1997), 1--39.

\bibitem{BLZm}
R.~{Bruggeman}, J.~{Lewis}, and D.~{Zagier}, \emph{{Period functions for Maass
  wave forms and cohomology}}, {Mem. Am. Math. Soc.} \textbf{1118} (2015),
  iii--v + 128.

\bibitem{BM09}
R.~{Bruggeman} and T.~{M\"uhlenbruch}, \emph{Eigenfunctions of transfer
  operators and cohomology}, Journal of Number Theory \textbf{129} (2009),
  158--181.

\bibitem{Chang_Mayer_transop}
C.-H. Chang and D.~Mayer, \emph{The transfer operator approach to {S}elberg's
  zeta function and modular and {M}aass wave forms for {${\rm PSL}(2,{\bf
  Z})$}}, Emerging applications of number theory ({M}inneapolis, {MN}, 1996),
  IMA Vol. Math. Appl., vol. 109, Springer, New York, 1999, pp.~73--141.

\bibitem{Choie_Zagier}
Y.~{Choie} and D.~{Zagier}, \emph{{Rational period functions for
  \(\rm{PSL}(2,\mathbb Z)\)}}, {A tribute to Emil Grosswald: number theory and
  related analysis}, Providence, RI: American Mathematical Society, 1993,
  pp.~89--108.

\bibitem{Deitmar_Hilgert}
A.~Deitmar and J.~Hilgert, \emph{A {L}ewis correspondence for submodular
  groups}, Forum Math. \textbf{19} (2007), no.~6, 1075--1099.

\bibitem{Folland}
G.~{Folland}, \emph{{Introduction to partial differential equations}}, 2nd ed.,
  Princeton, NJ: Princeton University Press, 1995.

\bibitem{Fraczek_Mayer_Muehlenbruch}
M.~Fraczek, D.~Mayer, and T.~M{\"u}hlenbruch, \emph{A realization of the
  {H}ecke algebra on the space of period functions for {$\Gamma_0(n)$}}, J.
  Reine Angew. Math. \textbf{603} (2007), 133--163.

\bibitem{Hassell}
A.~{Hassell}, \emph{{What is quantum unique ergodicity?}}, {Aust. Math. Soc.
  Gaz.} \textbf{38} (2011), no.~3, 158--167.

\bibitem{Hoer1}
L.~{H\"ormander}, \emph{{The analysis of linear partial differential operators.
  I: Distribution theory and Fourier analysis}}, reprint of the 2nd ed.,
  Berlin: Springer, 2003.

\bibitem{Katok_Ugarcovici}
S.~Katok and I.~Ugarcovici, \emph{Symbolic dynamics for the modular surface and
  beyond}, Bull. Amer. Math. Soc. (N.S.) \textbf{44} (2007), no.~1, 87--132.

\bibitem{Lewis}
J.~Lewis, \emph{Spaces of holomorphic functions equivalent to the even {M}aass
  cusp forms}, Invent. Math. \textbf{127} (1997), 271--306.

\bibitem{LZ_survey}
J.~{Lewis} and D.~{Zagier}, \emph{{Period functions and the Selberg zeta
  function for the modular group}}, {The mathematical beauty of physics: A
  memorial volume for Claude Itzykson. Conference, Saclay, France, June 5--7,
  1996}, Singapore: World Scientific, 1997, pp.~83--97.

\bibitem{LZ01}
\bysame, \emph{{Period functions for Maass wave forms. I}}, {Ann. Math. (2)}
  \textbf{153} (2001), no.~1, 191--258.

\bibitem{Mayer_thermo}
D.~Mayer, \emph{On the thermodynamic formalism for the {G}auss map}, Comm.
  Math. Phys. \textbf{130} (1990), no.~2, 311--333.

\bibitem{Mayer_thermoPSL}
\bysame, \emph{The thermodynamic formalism approach to {S}elberg's zeta
  function for {${\rm PSL}(2,{\bf Z})$}}, Bull. Amer. Math. Soc. (N.S.)
  \textbf{25} (1991), no.~1, 55--60.

\bibitem{Mayer_Muehlenbruch_Stroemberg}
D.~Mayer, T.~M{\"u}hlenbruch, and F.~Str{\"o}mberg, \emph{The transfer operator
  for the {H}ecke triangle groups}, Discrete Contin. Dyn. Syst. \textbf{32}
  (2012), no.~7, 2453--2484.

\bibitem{Moeller_Pohl}
M.~M{\"o}ller and A.~Pohl, \emph{Period functions for {H}ecke triangle groups,
  and the {S}elberg zeta function as a {F}redholm determinant}, Ergodic Theory
  Dynam. Systems \textbf{33} (2013), no.~1, 247--283.

\bibitem{Pohl_mcf_general}
A.~Pohl, \emph{A dynamical approach to {M}aass cusp forms}, J. Mod. Dyn.
  \textbf{6} (2012), no.~4, 563--596.

\bibitem{Pohl_mcf_Gamma0p}
\bysame, \emph{Period functions for {M}aass cusp forms for ${\Gamma}_0(p)$: A
  transfer operator approach}, Int. Math. Res. Not. \textbf{14} (2013),
  3250--3273.

\bibitem{Pohl_Symdyn2d}
\bysame, \emph{Symbolic dynamics for the geodesic flow on two-dimensional
  hyperbolic good orbifolds}, Discrete Contin. Dyn. Syst., Ser. A \textbf{34}
  (2014), no.~5, 2173--2241.

\bibitem{Ratcliffe}
J.~{Ratcliffe}, \emph{{Foundations of hyperbolic manifolds}}, 3rd expanded ed.,
  vol. 149, Cham: Springer, 2019.

\bibitem{Richards}
I.~Richards, \emph{Continued fractions without tears}, Math. Mag. \textbf{54}
  (1981), no.~4, 163--171.

\bibitem{Sarnak}
P.~{Sarnak}, \emph{{Recent progress on the quantum unique ergodicity
  conjecture}}, {Spectral geometry. Based on the international conference,
  Dartmouth, NH, USA, July 19--23, 2010}, Providence, RI: American Mathematical
  Society (AMS), 2012, pp.~211--228.

\bibitem{Selberg}
A.~Selberg, \emph{Harmonic analysis and discontinuous groups in weakly
  symmetric {R}iemannian spaces with applications to {D}irichlet series}, J.
  Indian Math. Soc. (N.S.) \textbf{20} (1956), 47--87.

\bibitem{Series}
C.~Series, \emph{The modular surface and continued fractions}, J. London Math.
  Soc. (2) \textbf{31} (1985), no.~1, 69--80.

\bibitem{Terras}
A.~{Terras}, \emph{{Harmonic analysis on symmetric spaces and applications.
  I}}, {New York etc.: Springer-Verlag}, 1985.

\bibitem{Venkov}
A.~{Venkov}, \emph{{Spectral theory of automorphic functions and its
  applications}}, Dordrecht etc.: Kluwer Academic Publishers, 1990.

\bibitem{Zelditch}
S.~{Zelditch}, \emph{{Recent developments in mathematical quantum chaos}},
  {Current developments in mathematics, 2009}, Somerville, MA: International
  Press, 2010, pp.~115--204.

\end{thebibliography}

\bibliographystyle{amsplain}

\setlength{\parindent}{0pt}

\end{document}